\documentclass [12pt]{article}
\usepackage[utf8]{inputenc}
\usepackage[english]{babel}
\usepackage{amsmath}
\usepackage[T1]{fontenc}
\usepackage{amsfonts, amssymb}
\usepackage[left=2.5cm,right=2.2cm,top=2.5cm,bottom=3.5cm]{geometry}
\usepackage{comment}
\usepackage{cmap}

\newtheorem{thm}{Theorem}
\newtheorem{lem}[thm]{Lemma}

\newtheorem{cor}[thm]{Corollary}
\newtheorem{example}[thm]{Example}
\newtheorem{rem}[thm]{Remark}

\newtheorem{cond}[thm]{Condition}

\def\R{{\mathbb R}}

\def\N{{\mathbb N}}

\DeclareMathOperator*{\esssup}{ess\,sup}

\DeclareMathOperator*{\argmax}{arg\,max}

\DeclareMathOperator*{\sign}{sign}

\DeclareMathOperator*{\Pe}{Pe}

\author{S.E. Nikitin}	
\title{Diffusions with resetting boundaries: ergodic control and renewal optimization}	
\date{\today}

\begin{document}

\maketitle

\begin{abstract}

We study ergodic optimization problems for one-dimensional diffusions with resetting boundaries, namely diffusions that are instantaneously reset according to prescribed distributions upon reaching designated boundary points. For controlled diffusions with resetting boundaries, we consider optimization of the ergodic cost over a broad class of admissible real-time controls. We derive a Hamilton--Jacobi--Bellman equation with non-local boundary conditions and establish a verification theorem for optimal stationary Markov controls. In the case of one resetting and one reflecting boundary, we also prove an existence result under weak integrability assumptions.
 
We then study renewal diffusions. Using their regenerative structure, we derive explicit formulas for the ergodic cost. These formulas are applied to optimization problems over families of drift coefficients, yielding explicit characterizations of optimal drift parameters in several examples. Finally, we study an infinite-penalty regime in which boundary hitting is prohibited and characterize the corresponding optimal ergodic cost.
\end{abstract}

\section{Introduction}
Optimization problems for controlled diffusions in both one-dimensional and higher-dimensional Euclidean spaces are well studied. The literature covers a wide range of models, including diffusions with reflecting boundaries, diffusions with Poisson jumps, switching diffusions, and many related settings. For an overview, see \cite{vB05}; for a more detailed treatment, see \cite{ABG12} and \cite{nK09}.

A further class of controlled diffusions, motivated by applications, incorporates resetting: at random times—typically occurring at a Poisson rate—the state is reset either to a fixed point or to a randomly chosen location. In contrast, optimization problems in which the resetting mechanism is triggered by hitting a boundary appear to have received much less attention. Such a problem was investigated in \cite{DRR21} for a finite interval in one dimension (including the case of a moving boundary) and for certain planar domains, using analytic methods. Similar models appear in financial mathematics, see \cite{DN90} and \cite{SS94}.

Related, though different in emphasis, is a substantial body of work on the ergodic behavior of diffusions that undergo random jumps away from the boundary. In particular, several papers establish exponential convergence to the stationary distribution for these processes; see \cite{GK02} for the one-dimensional case and \cite{BP09} for the multidimensional setting.  

In this work, we focus on ergodic optimization problems in one-dimensional diffusions with resetting boundaries. We study models in which the diffusion incurs a penalty and is reset at a random location upon hitting one of two boundary points, with the reset distribution allowed to depend on the boundary reached.

The paper is divided into two parts. The first part studies controlled diffusions with resetting boundaries and optimization of the ergodic cost over admissible stochastic controls that allow the diffusion to be steered in real time. The second part studies renewal diffusions, namely ordinary diffusions that are reset upon reaching a boundary, together with optimization problems over families of such diffusions.

The main contributions of this work are as follows. First, we derive a Hamilton--Jacobi--Bellman equation with non-local boundary conditions for controlled diffusions with resetting boundaries and characterize the optimal controls that maximize the ergodic cost via solutions to that equation. Second, in the case of one resetting and one reflecting boundary, we prove existence of solutions to the HJB equation under weak integrability assumptions, which yields existence of optimal stationary Markov controls. Third, we derive explicit formulas for the ergodic cost of renewal diffusions, and apply them to study optimization problems for renewal diffusions and analyze several representative examples. Finally, we consider the limiting infinite-penalty regime in which the diffusion is prohibited from reaching the boundary.

The paper is organized as follows.
In Section \ref{s:1dd}, we recall preliminary facts and notation concerning one-dimensional diffusions with reflecting boundaries.
In Section \ref{s:cdwrb}, we define controlled diffusions with resetting boundaries and formulate the optimization problem. 

In Section \ref{s:mr}, we define the Hamilton--Jacobi--Bellman equation with a non-local boundary condition and state the theorem characterizing optimal controls. In Section \ref{s:eofts}, we state an existence theorem. The results of these sections are proved in Section \ref{s:pcd}.

In Section \ref{s:rd}, we present explicit formulas for the ergodic cost of renewal diffusions. In Section \ref{s:prd}, we prove these results using the regenerative structure of renewal diffusions. 
In Section \ref{s:e}, we present examples of optimization over drift families. The proofs of these examples are given in Section \ref{s:pfoe}. 

In Section \ref{s:ip}, we study the infinite-penalty regime. The result of that section is proved in Section \ref{s:pi}.

\section{Preliminaries} \label{s:1dd}
We briefly recall the notion of one-dimensional diffusions with reflecting boundaries and fix notation used throughout the paper.
Let $W(\cdot )$ be a standard Brownian motion on $\R_{\ge 0}$. Let $-\infty \le A \le 0 < B$. 
 Let $K_{\mathrm{refl}} \subset \left\{ A,B \right\} \cap \R   $ be the set of reflecting boundary points and define 
 \[
	 \mathcal{X} := (A,B) \cup K_{\mathrm{refl}}.
 \]
 For convenience, if $A = -\infty$ then $[A,B] := (-\infty,B]$.
Let $\widetilde{\mu}: [A,B]\to \R$ be a drift coefficient  and  $\widetilde{\sigma}: [A,B]  \to \R_{>0}$ be a diffusion coefficient, both Borel measurable. 
A stochastic process $X(\cdot )$ is called a diffusion with reflecting boundaries in $K_{\mathrm{refl}}$, initial condition $x \in \mathcal{X}$, and killed upon leaving state space $\mathcal{X}$, if it satisfies the equation
\begin{eqnarray} \label{def:refl1}
\mathrm{d}X(t) = \widetilde{\mu} (X(t))\, \mathrm{d}t + \widetilde{\sigma} (X(t)) \,\mathrm{d}W(t) + \sum_{E \in K_{\mathrm{refl}}} r(E) \mathrm{d}\ell_{X}^{E}(t), 
\\ \nonumber  X(t) \in  \mathcal{X}, \quad t \in [0, T_1), \quad X(0) = x 
\end{eqnarray} 
where $r(A) := 1, \, r(B) := -1$, \, $T_1$ is the lifetime of $X(\cdot )$
\[
T_1 := \inf \left\{ t \ge 0 \, | \, X(t-) \in  \left\{ A,B \right\} \setminus K_{\mathrm{refl}}   \right\}, 
\]
and the reflection process $\ell_{X}^{E}(\cdot)$ is a continuous non-decreasing process that satisfies 
\[
\ell_{X}^{E}(t) = \int_{ 0 }^{t} 1_{\left\{ X(s) = E \right\} }  \, \mathrm{d}\ell_{X}^{E}(s), \quad t \ge 0.
\] 
By Tanaka's formula for semimartingales, the reflection process $\ell_{X}^{E}(\cdot)$ coincides with the local time at $E$ up to a conventional factor $\frac{1}{2}$.

For a definition of diffusions with reflecting boundaries in the sense of strong Markov processes with continuous paths (which are also called diffusions in the sense of Itô--McKean), see \cite[Ch.II.7]{BS02}.  

If $\widetilde{\mu}(\cdot )$ and $\widetilde{\sigma}(\cdot )$ satisfy the Engelbert--Schmidt condition 
\begin{equation} \label{cond:ES}
\frac{\widetilde{\mu}(\cdot )}{\widetilde{\sigma}(\cdot)^2}, \frac{1}{\widetilde{\sigma}(\cdot )^2}\in L_{1,\mathrm{loc}}[A,B] ,
\end{equation} 
then, for every initial condition $x \in \mathcal{X}$, the diffusion equation \eqref{def:refl1} admits a unique in law weak solution, that is, there is a solution on some probability space, and the law of that solution is unique, see \cite[Ref.~2.3]{UM12} and the references there.  

\section{Controlled diffusions with resetting boundaries} \label{s:cdwrb}
Let $\left( \Omega, \mathcal{F}, \mathbb{P} \right)$ be a probability space equipped with a filtration $ \left( \mathcal{F}_{t} \right)_{t \ge 0}$, and let $W(\cdot )$ be a Brownian motion adapted to this filtration. 
 Let  $U$ be a compact metric space, and let $ \left\{ u(t)\right\}_{t \ge 0} $ be an $ \left( \mathcal{F}_{t} \right)$-adapted $U$-valued process, called a control process. A tuple
\[
	\left(  \Omega, \mathcal{F}, \mathbb{P}, \left( \mathcal{F}_{t} \right)_{t \ge 0}, W(\cdot ), u(\cdot ) \right)
\]
is called a {\it weakly admissible control}. This corresponds to a weak formulation, since the probability space is part of the control.
Let $K_{\mathrm{refl}} \subset \left\{ A,B \right\}  \cap \R$ be the set of reflecting boundary points as before, and let $K_{\mathrm{res}} := \left( \left\{ A,B \right\} \cap \R  \right) \setminus  K_{\mathrm{refl}} $ be the set of resetting boundary points. We assume that $K_{\mathrm{res}} \neq \varnothing$. 
Let $\left\{ \nu_{E} \, | \, E \in K_{\mathrm{res}}\right\}$ be probability distributions on $\mathcal{X}$, and let
\[
	\left\{ \gamma_n(E) \sim \nu_{E} \, | \, n \in \N_{> 0}, E \in K_{\mathrm{res}}  \right\}
\]
be a collection of independent random variables. The reset variables are assumed to be independent of $W(\cdot )$ and of the pre-reset filtration at the corresponding reset times. At a reset time, when the process hits a boundary point $E \in K_{\mathrm{res}}$, it is reset to a new point drawn according to $\nu_E$, represented by the random variable $\gamma_n(E)$.

Let $\mu: [A,B] \times U \to \R$ and $\sigma: [A,B] \times U \to \R_{> 0}$ be Borel measurable functions.
For a given weakly admissible control, if there exists an $\left( \mathcal{F}_{t} \right)$-adapted pathwise (strong) solution to the following piecewise-defined stochastic differential equation, we call that solution the {\it controlled diffusion with resetting boundaries} in $K_{\mathrm{res}}$:
\begin{eqnarray} \nonumber
\mathrm{d}Z_{u}(t) = \mu(Z_{u}(t),u(t)) \mathrm{d}t + \sigma(Z_{u}(t),u(t))\mathrm{d}W(t)  + \sum_{E \in K_{\mathrm{refl}}} r(E) \mathrm{d}\ell_{Z_{u}}^{E}(t), \\ \label{def:Z_u}
\quad t \in [T_{u}(n-1), T_{u}(n)), \quad
Z_{u}(T_{u}(n)) = \gamma_{n}(E_{u}(n)), \quad n \in \N_{> 0}, \quad Z_{u}(0) = 0,   
\end{eqnarray} 
where $T_{u}(n)$ denotes the $n$-th time at which $Z_{u}(\cdot )$ hits $K_{\mathrm{res}}$
\[
	T_{u}(n) := \inf \left\{ t > T_{u}(n-1) \, | \, Z_{u}(t-) \in K_{\mathrm{res}}    \right\}, \quad n \in \N_{> 0}, \quad T_{u}(0) := 0,  
\] 
and $E_{u}(n)$ denotes the boundary point that is hit at time $T_{u}(n)$ before reset 
\[
E_{u}(n) := Z_{u}(T_{u}(n)-), \quad n \in \N_{> 0}. 
\] 
The process $Z_{u}(\cdot )$ is defined on $[0, T_{u}(\infty))$, where $T_{u}(\infty) := \lim\limits_{n \to \infty} T_{u}(n)$, and
\[
Z_{u}(t) \in \mathcal{X}, \quad t \in [0, T_{u}(\infty)).
\]
\subsection{The optimization problem} \label{s:tof}
Let $\mathcal{D}^{\circ}$ denote the class of all weakly admissible controls for which equation \eqref{def:Z_u} admits a pathwise (strong) solution. Let $\mathcal{D}$ be the subclass of controls in $\mathcal{D}^{\circ}$ satisfying $T_{u}(\infty) \stackrel{\text{a.s.}}{=} \infty$. 
For brevity, we write $u \in \mathcal D$ (or $\mathcal{D}^{\circ} $) whenever
\[
\left(\Omega,\mathcal F,\mathbb P,(\mathcal F_t)_{t\ge0},W(\cdot),u(\cdot)\right)
\in\mathcal D \quad (\text{or }\mathcal{D}^{\circ} ).
\]

Let $u \in \mathcal{D}$.
Define $N_u(T)$ as the number of resets up to and including time $T \ge 0$, namely
\[
N_u(T) := \sup \left\{ n \in \N_{\ge 0} \, | \, T_{u}(n) \le T \right\}.  
\]
For $E \in K_{\mathrm{res}} $, define $q_{u}(E,T)$ as the number of resets caused by reaching the boundary point $E$ up to and including time $T \ge 0$, namely 
\begin{eqnarray*}
	q_{u}(E,T) &:=& \#\left\{ i \in \left\{ 1, \ldots, N_{u}(T) \right\}  \, | \, E_{u}(i) = E  \right\}.
\end{eqnarray*} 

 Let $f: [A,B]\times U \to \R$ be a Borel measurable function, called the reward function. Define the pathwise long-run average (ergodic) {\it objective functional} with penalty vector $P_{K_{\mathrm{res}}} := \left\{ P_{E} \in \R \, | \, E \in K_{\mathrm{res}}   \right\} $, where a penalty $P_{E}$ is incurred when the process hits $E \in K_{\mathrm{res}}$, as
\[ 
	J(u,f, P_{K_{\mathrm{res}}}) :=   \limsup_{T \to \infty} \frac{1}{T} \left( \int_{ 0 }^{ T }  f(Z_{u}(t),u(t)) \mathrm{d}t - \sum_{E \in K_{\mathrm{res}} } P_{E}q_{u}(E, T) \right). 
\] 
The optimization problem considered in this paper is to maximize the objective functional
\[
J(u,f,P_{K_{\mathrm{res}}})
\]
{\it almost surely} over all admissible controls \(u\in\mathcal D\).
\begin{rem}
One can also consider an optimization problem based on a deterministic version of the objective functional $J$, obtained by replacing the pathwise $\limsup$ with the $\limsup$ of expectations. Analogous results hold for this deterministic formulation. They are weaker than the corresponding pathwise results, but their proofs require less technical work.
\end{rem}

\section{Main results for controlled diffusions with resetting boundaries on a bounded interval} \label{s:mr}
\subsection{Standing assumptions} \label{s:sa}
In this section we restrict attention to the case $A > -\infty$.
The next condition provides integrability and growth assumptions on the drift, diffusion, and reward coefficients. These assumptions will be used both in the analysis of the HJB equation and in the construction of optimal controls.
\begin{cond} \label{cond:MSF_full}

The functions $\mu(x,\eta)$, $\sigma(x,\eta)$, and $f(x,\eta)$ are continuous in $\eta$, measurable in $x$, and there exist positive Borel measurable functions $\Lambda, \varsigma, F$ on $[A,B]$ such that, for a.e. $x \in (A,B)$, for all $\eta \in U$, it holds that 
\begin{eqnarray*}
	&&\Lambda(x) \ge \frac{|\mu(x,\eta)|}{\sigma^2(x,\eta)}, \quad	F(x) \ge \frac{|f(x,\eta)|}{\sigma^2(x,\eta)}, \\
	&& \ \varsigma(x)\le \sigma^2(x,\eta), 
\end{eqnarray*}
and
\begin{equation*} \label{cond:MSF}
\Lambda, F, \frac{1}{\varsigma} \in L_{1}[A,B].
\end{equation*} 
\end{cond}
Condition \ref{cond:MSF_full} may be viewed as a controlled analogue of the Engelbert--Schmidt condition \eqref{cond:ES}.

The next lemma shows that, on a bounded interval, Condition \ref{cond:MSF_full} guarantees non-explosion of the reset times.
In particular, it implies that the controlled diffusions corresponding to controls in $\mathcal{D}^{\circ}$ are well defined for all times almost surely.
\begin{lem} \label{l:T}
 Let $u \in \mathcal{D}^{\circ}$. 
If Condition \ref{cond:MSF_full} holds, then there exists $c > 0$ such that, with probability 1, for all sufficiently large $n \in \N_{>0}$, 
\[
T_{u}(n) \ge cn.
\]
In particular, we have $\mathcal{D}^{\circ} = \mathcal{D}$.
\end{lem}

\subsection{Optimal controls and the Hamilton--Jacobi--Bellman equation with non-local boundary conditions} \label{s:hjb}

For $\alpha \in \R$, consider the Hamilton--Jacobi--Bellman (HJB) equation with non-local boundary conditions, for a.e. $x \in (A,B)$,
\begin{equation} \label{e:HJB}
	\begin{cases} 
		\alpha = \sup\limits_{\eta \in U } \left( \frac{1}{2} \sigma^2(x,\eta) V_{\alpha}^{\prime \prime}(x) + \mu(x,\eta) V^{\prime} _{\alpha}(x) + f(x,\eta)\right) , 
		\\ 
		V_{\alpha}(E) + P_{E} =\int_{\mathcal{X}} V_{\alpha}(x) \mathrm{d}\nu_{E}(x), \quad E  \in K_{\mathrm{res}}, 
		\\
		V_{\alpha}^{\prime}(E) = 0, \quad E \in K_{\mathrm{refl}}.
	\end{cases}
\end{equation}

A control $u$ is called a {\it stationary Markov (feedback) control}
if there exists a Borel measurable function $h:[A,B]\to U$ such that
\[
u(t)=h(Z_u(t)), \quad \text{for all } t \ge 0, \text{ a.s.}
\]
The following theorem shows that any sufficiently regular solution to the HJB equation \eqref{e:HJB} yields an optimal stationary Markov control.
\begin{thm} \label{t:main}
Assume
that Condition \ref{cond:MSF_full} holds. Suppose moreover that there exists $\alpha \in \R$ such that equation \eqref{e:HJB} admits a solution $V_{\alpha} \in W^{2,1}[A,B]$ (equivalently, $V_{\alpha}^{\prime}\in AC[A,B]$).
If $u \in \mathcal{D}$,
then
\[
J(u,f,P_{K_{\mathrm{res}}})  \stackrel{\text{a.s.}}{\le} \alpha.
\]
There exists a measurable selector $\eta^{*}: [A,B] \to U$ satisfying 
\begin{equation} \label{def:u0*}
\eta^{*}(x) \in \argmax\limits_{\eta \in U} \left( \frac{1}{2} \sigma^2(x,\eta) V_{\alpha}^{\prime \prime}(x) + \mu(x,\eta) V_{\alpha}^{\prime}(x) + f(x,\eta)\right), \quad \text{for a.e. }x  \in (A,B).
\end{equation} 
There exist a filtered probability space and an adapted process $Z_{u^{*}}(\cdot )$, driven by a Brownian motion on that space, satisfying \eqref{def:Z_u} under the stationary Markov  control
\begin{equation} \label{def:u*}
	u^{*}(t) := \eta^{*}(Z_{u^{*}}(t)), \quad t \ge 0
\end{equation}
such that 
\begin{equation} \label{e:main_equality}
	J(u^{*}, f, P_{K_{\mathrm{res}}}) \stackrel{\text{a.s.}}{=}  \lim_{T \to \infty} \frac{1}{T} \left( \int_{ 0 }^{ T }  f(Z_{u^{*}}(t),u^{*}(t)) \mathrm{d}t - \sum_{E \in K_{\mathrm{res}} } P_{E}q_{u^{*}}(E,T) \right) \stackrel{\text{a.s.}}{=}\alpha.
\end{equation} 
\end{thm}
\begin{rem}
 One may also consider controlled diffusions with partially reflecting and partially resetting boundary points. In that case, for the ergodic cost, one obtains a similar HJB equation with a non-local boundary condition involving the derivative at the boundary.  
\end{rem}

\subsection{Canonical drift-control problem}
A natural example is the direct control of the drift coefficient, constrained to take values in a bounded interval.
\begin{example} \label{exa:0}
Let $U := [-1,1]$, $\mu(x,\eta) := \eta$, and $\sigma(x,\eta) := 1$ for $(x,\eta) \in [A,B] \times U$. Let $f \in L_1[A,B]$ and $P_{E} \in \R$. 
Then, for $\alpha \in \R$, the HJB equation \eqref{e:HJB} takes the form, for a.e. $x \in (A, B)$, 
\begin{eqnarray*}
\begin{cases}
	-V^{\prime \prime}(x) = 2 \left( |V^{\prime }(x)| + f(x) - \alpha \right) ,
\\
	V_{\alpha}(E) + P_{E} =\int_{\mathcal{X}} V_{\alpha}(x) \mathrm{d}\nu_{E}(x), \quad E  \in K_{\mathrm{res}}, 
		\\
		V_{\alpha}^{\prime}(E) = 0, \quad E \in K_{\mathrm{refl}}.
\end{cases}
\end{eqnarray*} 
and the optimal drift coefficient is given by 
\[
	\mu^{*}(x) = \sign(V^{\prime}(x)), \quad  x \in [A,B]. 
\] 
\end{example}
Thus the optimal drift coefficient is of bang-bang type: at each state $x \in [A,B]$, the controller applies the maximal admissible drift in the direction determined by the sign of $V^{\prime} (x)$.

\subsection{Existence of solutions to the HJB equation in the one-reflecting/one-resetting case} \label{s:eofts}
We next show that, in the case of one reflecting and one resetting boundary, the HJB equation \eqref{e:HJB} admits a solution.
\begin{thm} \label{t:equ}
	Under Condition \ref{cond:MSF_full}, suppose that $K_{\mathrm{res}} = \left\{ B  \right\} $ and $K_{\mathrm{refl}} = \left\{ A \right\} $.

Then there exist $\alpha \in \R$ and a function $V_{\alpha}\in W^{2,1}[A,B]$  such that $V_{\alpha}$ is a solution to \eqref{e:HJB}. 

\end{thm}

Theorem \ref{t:equ} establishes existence of solutions to the HJB equation under Condition \ref{cond:MSF_full}. Consequently, Theorem \ref{t:main} applies in the one-reflecting/one-resetting case.
Condition \ref{cond:MSF_full} is natural in our framework since, together with the Engelbert--Schmidt condition \eqref{cond:ES}, it ensures that a solution to the HJB equation determines an admissible control in $\mathcal D$. Extending the existence result to the case of two resetting boundaries appears to require substantially different techniques, such as fixed-point methods under more standard, but stronger, regularity assumptions, and is therefore left for future work.

Finally, we state an auxiliary existence result needed in the proof of Theorem \ref{t:main}. The following lemma provides a Lyapunov-type function used to control the local martingale terms appearing in the Itô decomposition.
\begin{lem} \label{l:equ}
Under Condition \ref{cond:MSF_full}, there exist $\Gamma \in \R$ and a function $H \in W^{2,1}[A,B]$ such that, for a.e. $x \in (A,B)$,
\begin{equation}
	\begin{cases} \label{e:H}
		\Gamma \ge \sup\limits_{\eta \in U } \left( \frac{1}{2} \sigma^2(x,\eta) H^{\prime \prime}(x) + \mu(x,\eta) H^{\prime}(x) + \sigma^2(x,\eta)\right) , 
		\\ 
		H(E) \ge \int_{\mathcal{X}} H(x) \mathrm{d}\nu_{E}(x), \quad E  \in K_{\mathrm{res}}, 
		\\
		H^{\prime}(E)r(E) \le 0, \quad E \in K_{\mathrm{refl}}.
	\end{cases}
\end{equation} 
\end{lem}

\section{Main results for renewal diffusions} \label{s:rd}
Renewal diffusions form a natural class of stochastic processes. They arise by repeatedly resetting a diffusion upon reaching designated boundary points, thereby generating a regenerative structure. Such processes appear in a variety of settings involving restart mechanisms, and their long-run behavior can often be analyzed explicitly using renewal-theoretic methods.

In this paper, renewal diffusions also arise naturally from the control problem studied in the previous sections. Indeed, Theorem \ref{t:main} shows that, in many cases, an optimal control is a stationary Markov control. Under a fixed stationary Markov control, the controlled diffusion becomes a renewal diffusion. This connection motivates the study of renewal diffusions from the perspective of ergodic optimization.

The regenerative structure of renewal diffusions makes it possible to derive explicit formulas for the ergodic cost.  We will subsequently use these formulas to study optimization problems over families of drift coefficients.

Let $\mu: [A,B] \to \R$ and $\sigma: [A,B] \to \R_{>0}$ be Borel measurable functions satisfying Engelbert--Schmidt condition \eqref{cond:ES}, and let $f: [A,B] \to \R$ be a Borel measurable function. Let 
$K_{\mathrm{res}}$ be either $\left\{ A,B \right\} $ or $\left\{ B \right\} $.

Let the control space be trivial, that is, $U := \left\{ 0 \right\} $. Since there is no dependence on $u$, we omit it from the notation. Let $Z(\cdot )$ be the controlled diffusion with resetting boundaries satisfying \eqref{def:Z_u}. 
In this case, we call  $Z(\cdot )$ a {\it renewal diffusion}. In the single-boundary
case, after the first reset, the successive excursions are i.i.d.; in the two-boundary
case, the successive excursions are governed by the embedded boundary-hit chain. Unlike in the controlled setting, non-explosion of reset times for renewal diffusions follows directly from this regenerative structure.
\begin{lem} \label{l:d}
$T(\infty) \stackrel{\text{a.s.}}{=} \infty $. 
\end{lem}

We consider the pathwise long-run average (ergodic) objective functional with penalty $P_{K_{\mathrm{res}}} := \left\{ P_{E} \in \R \, | \, E \in K_{\mathrm{res}}   \right\} $ upon reaching the boundary point $E \in K_{\mathrm{res}}$ defined by 
\begin{equation*} \label{def:JK}
	J(f,P_{K_{\mathrm{res}}}) :=   \lim_{T \to \infty} \frac{1}{T} \left( \int_{ 0 }^{ T }  f(Z(t)) \mathrm{d}t - \sum_{E \in K_{\mathrm{res}} } P_{E}q(E,T) \right). 
\end{equation*} 
We use $\lim$ instead of $\limsup$ as in Section \ref{s:tof}. Such a limit typically exists for renewal diffusions, in analogy with the optimal stationary Markov control in Theorem  \ref{t:main}.   

\subsection{Green function representation} \label{s:gf}

The following Green function representation allows explicit computation of the ergodic objective functional.
The {\it speed measure} of a diffusion $X(\cdot )$ has density 
\begin{equation*} \label{def:m_full}
   m(x) := \frac{1}{\sigma^2(x)}\exp\left(2 \int_{0}^{x} \frac{\mu(s)}{\sigma^2(s)} \, \mathrm{d}s\right), \quad x \in [A,B].
\end{equation*} 
The {\it scale function} $S(\cdot )$ of $X(\cdot )$ is defined as
\begin{equation*} \label{def:S}
	S(x) :=  \int_{0}^{x} \exp \left( -2 \int_{0}^{y} \frac{\mu(s)}{\sigma^2(s)} \mathrm{d}s \right) \,\mathrm{d}y, \quad x \in [A,B].
\end{equation*} 
These functions are fundamental characteristics of one-dimensional diffusions, see \cite[Ch.~IV.11, Ch.~IV.15]{aB17} and \cite[Ch.~5.5.C]{KS98}.

Define the Green functions for a diffusion killed upon reaching $K_{\mathrm{res}}$ as 
\begin{eqnarray*}
	G_{AB}(x,y) &:=&  \frac{\left( S(x) - S(A) \right)\left( S(B)-S(y)  \right) }{S(B)-S(A)}, \quad A \le x \le y \le B, \quad K_{\mathrm{res}} = \left\{ A,B \right\}, 
\\
	G_{B}(x,y) &:=& S(B)-S(y), \quad A \le x \le y \le B, \quad K_{\mathrm{res}} = \left\{ B \right\}, 
\end{eqnarray*} 
and extend them symmetrically to the region $x > y$. We assume that all values of the scale function appearing in these formulas are finite.
Define the associated Green functional as
\[
G_{K_{\mathrm{res}}}(x,f) := 2\int_{A}^{B} G_{K_{\mathrm{res}}}(x,y) f(y) m(y) \, \mathrm{d}y.
\] 
Define the averaged Green functional as
\begin{equation} \label{def:GKE}
	G_{K_{\mathrm{res}}}^{E}(f) :=  \int_{A}^{B} G_{K_{\mathrm{res}}}(x,f) \, \mathrm{d}\nu_{E}(x), \quad E \in K_{\mathrm{res}}. 
\end{equation} 
These objects naturally arise in the representation of long-run average rewards for regenerative diffusions.
For brevity, we write the subscripts $AB$ and $B$ in place of $K_{\mathrm{res}} = \left\{ A,B \right\} $ and $K_{\mathrm{res}} = \left\{ B \right\} $, respectively.
We also write $G_{B}(f)$ instead of $G_{B}^{B}\left(f \right) $ and $J(f,P_{B})$ instead of $J(f,P_{K_{\mathrm{res}}})$. Since the following results are stated for fixed $\mu(\cdot )$ and $\sigma(\cdot )$, we omit them from the notation.  
\subsection{Explicit formulas for the ergodic cost} \label{s:ef}
We first consider the case of two resetting boundaries on a bounded interval. In this setting, the long-run average objective functional can be expressed explicitly in terms of the averaged Green function and the asymptotic frequencies of hitting the two boundaries.
\begin{thm} \label{t:AB}
	Suppose that $A > -\infty$. Let $K_{\mathrm{res}} := \left\{ A,B \right\} $, and assume that $|S(E)| < \infty $ for $ E \in \left\{ A,B \right\} $.
Suppose that $\frac{f}{\sigma^2} \in L_1[A,B]$. 
Then
\begin{equation} \label{e:JfAB}
J(f,P_{AB}) \stackrel{\text{a.s.}}{=} \frac{\pi_{A} \left( G_{AB}^{A}(f) -P_{A}\right)  + \pi_{B} \left( G_{AB}^{B}(f ) - P_{B} \right) }{\pi_{A} G_{AB}^{A}(1) + \pi_{B} \, G_{AB}^{B}(1)} =: \alpha,  
\end{equation} 
where
\[
\pi_{A} := \frac{p_{BA}}{p_{BA}+p_{AB}}, \quad \pi_{B} := \frac{p_{AB}}{p_{BA}+p_{AB}}, 
\] 
\begin{eqnarray*}
	p_{BA} := \int_{A}^{B} \frac{ S(B) - S(x) }{S(B)-S(A)}   \, \mathrm{d}\nu_{B}(x), \quad 
	p_{AB} := \int_{A}^{B} \frac{ S(x) - S(A) }{S(B)-S(A)}   \, \mathrm{d}\nu_{A}(x). 
\end{eqnarray*} 
\end{thm}
\begin{rem}
When $K_{\mathrm{res}} = \left\{ A,B \right\}$, the empirical frequencies of hitting $A$ and $B$ converge almost surely to $\pi_{A}$ and $\pi_{B}$, respectively. 
\end{rem}

We next consider the case of a single resetting boundary, with the other boundary either reflecting or located at $-\infty$. In contrast to the case of two resetting boundaries, the recurrence properties of the underlying diffusion now play a crucial role.
\begin{thm} \label{t:B}
	Let $K_{\mathrm{res}}:={\left\{B  \right\} }$. 
	Suppose that $f \ge 0$,\ $|S(x)| < \infty$ for all $x \in [A,B]$, and, if $A = -\infty$, then $S(A) = -\infty$.
If either  $G_{B}(f) < \infty$ or $G_{B}(1) < \infty$, then 
\[
J(f,P_{B}) \stackrel{\text{a.s.}}{=} \frac{G_{B}(f)-P_{B}}{G_{B}(1)}.
\] 
Here we adopt the conventions $\frac{c_1}{\infty} = 0, \frac{\infty}{c_2} = \infty,$ for $c_1 \in \R, c_2 > 0$.
\end{thm}
\subsection{Exceptional cases for a single resetting boundary}
The following theorem describes exceptional cases in which the ergodic objective functional becomes infinite or degenerates to the limiting value of the reward function at $-\infty$.
We require the following two conditions.
\begin{eqnarray} \label{cond:fs}
&&\frac{f(\cdot )}{\sigma^2(\cdot )} \not \in L_{1,\mathrm{loc}}[A,B], 
\\ \label{cond:fS} 
	\forall \varepsilon > 0 \, && \int_{B-\varepsilon}^{B} (S(B)-S(x)) f(x) m(x) \, \mathrm{d}x = \infty. 
\end{eqnarray}
\begin{thm} \label{t:special}
	Let  $K_{\mathrm{res}}:=\left\{ B \right\} $. 
	Suppose that $f \ge 0$, \ $S(x) < \infty$ for all $x \in (A,B]$ and, if $A = -\infty$, then $S(A) = -\infty$. 
If either condition \eqref{cond:fs} or \eqref{cond:fS} holds, then  
\[
J(f,P_{B} ) \stackrel{\text{a.s.}}{=} \infty.
\] 
If $G_B(1) =  \infty$, $f(\cdot )$ has a limit at $-\infty$, and neither condition \eqref{cond:fs} nor \eqref{cond:fS} holds, then
\[
J(f,P_{B} ) \stackrel{\text{a.s.}}{=} \lim_{x \to -\infty} f(x).
\] 
\end{thm}

\section{Optimization over drift families} \label{s:e}
We now study optimization problems over families of drift coefficients.
We consider only the case $\nu_{A}=\nu_{B}=\delta_{0}$. Throughout this section the diffusion coefficient is fixed as $\sigma \equiv 1$, and optimization is performed only over the drift coefficient.
Let
\[
	\mathcal{D}^{\prime}  \subset \left\{\mu \in L_{1, \mathrm{loc}}[A,B]\right\}
\]
be a set of drift coefficients.
Each  $ \mu \in \mathcal{D}^{\prime}$ determines a renewal diffusion with drift coefficient $\mu$ and deterministic reset to 0 upon reaching $K_{\mathrm{res}} \subset \left\{ A,B \right\} $, where $K_{\mathrm{res}}$ is fixed throughout.

The drift-control problem of Example \ref{exa:0} may be viewed as an optimization problem over the family 
 \[
	 \mathcal{D}_{0} = \left\{ \mu: [A,B] \to [-1,1] \right\}. 
\] 
There we showed that the optimal drift is of bang-bang type and that the corresponding stationary Markov control is optimal over a much wider class of stochastic controls.

From now on, we fix $K_{\mathrm{res}} := \left\{ B \right\} $, $P_{B} \ge 0$, and $f(x) := \frac{1}{B} x \cdot  1_{[0, B]}(x)$. We write
\[
J(\mu, f, P_{B}) := J(f,P_B), \quad  G_{B}(\mu, f) := G_{B}^{B}(f),
\] 
emphasizing the dependence on the drift coefficient.
We are interested in the maximal value of the objective functional over renewal diffusions in $\mathcal{D}^{\prime} $, namely 
\[
I(\mathcal{D}^{\prime} , f, P_{B}) := \sup_{ \mu \in \mathcal{D}^{\prime} } \mathbb{E} \, J(\mu, f,P_{B}). 
\] 
For all drift coefficients $\mu$ considered below, the quantity $J(\mu, f,P_{B})$ is almost surely constant. We therefore take expectation only to avoid dealing with null sets when taking the supremum over an uncountable set.

We now consider several simple families of drift coefficients for which the optimal value can be characterized explicitly. Unlike Example \ref{exa:0}, the following examples can be treated using the explicit formulas derived in Section \ref{s:rd}.

The first example illustrates the trade-off between rapid resetting and the penalty incurred at the resetting boundary. For small penalties, positive drift increases the frequency of resetting and improves the average reward. When the penalty exceeds the critical value $\frac{B^2}{3}$, resetting becomes too costly and the optimal strategy is to suppress boundary hitting by choosing zero drift.
\begin{example} \label{exa:1}

Let $A = -\infty$. Consider the set $\mathcal{D}_{1} =  \left\{ \mu \equiv \mathrm{const} \ge 0\right\} $. Each $\mu \in \mathcal{D}_{1}$ corresponds to a Brownian motion with a constant drift $\mu \ge 0$, a single resetting boundary at $B$, and deterministic reset to $0$ upon reaching $B$.
Then, for $P_{B} = 0$, we have 
\begin{equation} \label{case:0} 
I(\mathcal{D}_{1},f,P_{B}) = \lim_{\mu \to \infty} \mathbb{E} J(\mu,f, 0) = \frac{1}{2}.
\end{equation} 

For $P_{B} \in (0, \frac{B^2}{3})$, we have 
\begin{equation*} \label{case:1}
I(\mathcal{D}_{1},f,P_{B}) = \frac{2\Pe^2 - 2\Pe + 1 - e^{-2\Pe}}{4\Pe^2} - \frac{P_{B} \Pe}{B^2},
\end{equation*} 
where the P\'eclet number $\Pe > 0$ is the unique solution to the equation
\begin{equation} \label{case:2}
\frac{\Pe (1+e^{-2\Pe}) + e^{-2\Pe} - 1}{2\Pe^3} = \frac{P_{B}}{B^2}.
\end{equation} 
In that case, $I(\mathcal{D}_{1},f,P_{B}) = \mathbb{E} J(\mu,f, P_{B})$, where $\mu = \frac{\Pe}{B}$. 

For $P_{B} \ge \frac{B^2}{3}$, we have 
\[
	I(\mathcal{D}_1,f,P_{B}) = \mathbb{E} J(0, f, P_{B}) = 0.
\]
\end{example}

The next example was considered in  \cite{DRR21}. Compared with Example \ref{exa:1}, the reflecting boundary at $0$ prevents the diffusion from escaping to $-\infty$, so both positive and negative constant drifts can be considered. The optimal drift balances spending time near the upper boundary, where the reward is larger, against the penalty incurred upon resetting.
\begin{example} \label{exa:2}
Let $A = 0$. 
Consider the set $\mathcal{D}_{2} :=  \left\{ \mu \equiv \mathrm{const} \in \R\right\} $. Each $\mu \in \mathcal{D}_{2}$ corresponds to a Brownian motion with a constant drift $\mu \in \R$, a resetting boundary at $B$, a reflecting boundary at $0$, and deterministic reset to $0$ upon reaching $B$.
We have 
\begin{equation} \label{e:IfPe}
I(\mathcal{D}_{2},f,P_{B}) = \sup_{\Pe \in \R} \left\{ \frac{(2\Pe^2 - 2\Pe + 1)e^{2\Pe} - 1}{2\Pe((2\Pe - 1)e^{2\Pe} + 1)} - \frac{P_{B}}{B^2} \cdot  
\frac{2\Pe^2}{2\Pe - 1 + e^{-2\Pe}} \right\} .
\end{equation} 
The optimal drift is given by $\mu = \frac{\Pe}{B}$, where $\Pe \in \R$ maximizes \eqref{e:IfPe}. 
\end{example}
\section{Infinite-penalty regime} \label{s:ip}
This limiting problem complements the finite-penalty resetting model by describing what happens when boundary hitting is made prohibitively costly. We consider the case in which hitting the boundary is forbidden, corresponding formally to the regime $P_{A},P_{B} \to \infty$ in the resetting model. Optimization problems similar to those considered in Section \ref{s:e} were studied in \cite{LS15} and \cite{LN24} in the context of energy-conserving approximations of Brownian motion, constrained to remain pathwise close to, or above, the Brownian motion.

Let $f: (A, B) \to \R_{\ge 0}$ be a Borel measurable function. Let $X(\cdot )$ be a diffusion with state space $(A,B)$. Consider the objective functional 
\[
    J(X,f) := \mathbb{E} \left( \lim_{T \to \infty} \frac{1}{T} \int_{0}^{T} f(X(t)) \, \mathrm{d}t \right) .  
\] 
In all cases considered, the limit exists and is almost surely constant.

Denote by $\mathcal{D}_{AB}$ the class of homogeneous diffusions, started at 0,  with {\it continuous} on $(A,B)$ drift coefficient $\mu(\cdot )$ and constant diffusion coefficient $\sigma \equiv 1$, that almost surely do not leave the interval $(A, B)$, and for which the above objective functional is well-defined (possibly infinite).
We are interested in the maximal value of the objective functional $J(X,f)$ over $X \in \mathcal{D}_{AB}$. Define 
\[
I(\mathcal{D}_{AB},f) := \sup \left\{ J(X,f) \, | \, X \in \mathcal{D}_{AB}\right\}.  
\] 

Equivalently, this problem can be formulated as follows: the diffusion earns a reward of $\int_{T_1}^{T_2} f(X(t)) \, \mathrm{d}t $ over the interval $[T_1, T_2]$ provided it does not reach levels $A$ and $B$; upon reaching either $A$ or $B$, it incurs a penalty of $-\infty$.  

\begin{thm} \label{t:infinite} $I(\mathcal{D}_{AB},f) = \esssup f.$ \end{thm}
\begin{rem}
	The quantity $\esssup f$ is the largest possible value of the objective functional. The theorem shows that this universal upper bound can be approached arbitrarily closely by diffusions in the class $\mathcal{D}_{AB}$.
\end{rem}
\begin{rem}
The supremum in Theorem \ref{t:infinite} is generally not attained by any diffusion. 
\end{rem}

\section{Proofs of the results for controlled diffusions with resetting boundaries} \label{s:pcd}
\subsection{Proof of Theorem \ref{t:main}}
{\bf Itô decomposition}.
Let $u \in \mathcal{D}$ and define 
\[
\mathcal{L}_{u}V_{\alpha}(t) := \frac{1}{2}\sigma^2(Z_{u}(t), u(t)) V_{\alpha}^{\prime \prime} (Z_{u}(t)) + \mu(Z_{u}(t),u(t))V^{\prime} _{\alpha}(Z_{u}(t)), \quad t \ge 0.
\] 
For $T \ge 0$, by using \eqref{e:HJB} and generalized Itô's formula for $V_{\alpha} \in W^{2,1}[A,B]$, we obtain 
\begin{eqnarray} \nonumber
&& \int_{ 0 }^{T}  f(Z_{u}(t),u(t)) \, \mathrm{d}t - \sum_{E \in K_{\mathrm{res}} } P_{E}q_{u}(E, T) 
	\\ \label{e:Ito_inequality} &\le& \int_{ T_{u}(N_u(T))}^{T} \left( \alpha - \mathcal{L}_{u}V_{\alpha}(t) \right)  \, \mathrm{d}t + \sum_{i=1}^{N_u(T)}  \int_{ T_{u}(i-1) }^{T_{u}(i)} \left( \alpha - \mathcal{L}_{u}V_{\alpha}(t) \right)  \, \mathrm{d}t - \sum_{E \in K_{\mathrm{res}} } P_{E}q_{u}(E, T) \qquad 
	\\ \nonumber&=&
	\alpha T  - V_{\alpha}(Z_{u}(T)) + V_{\alpha}(0)  + \int_{ T_{u}(N_u(T)) }^{T} V^{\prime} _{\alpha}(Z_{u}(t)) \, \sigma(Z_{u}(t),u(t))  \, \mathrm{d}W(t) 
	\\ \nonumber &+& \sum_{i=1}^{N_u(T)} \left( V_{\alpha}(Z_{u}(T_{u}(i))) - V_{\alpha}(E_{u}(i))  - P_{E_{u}(i)}  + \int_{ T_{u}(i-1) }^{T_{u}(i)} V^{\prime} _{\alpha}(Z_{u}(t)) \, \sigma(Z_{u}(t),u(t))  \, \mathrm{d}W(t)  \right)	
	\\ \label{e:Ito} &+& \sum_{E \in K_{\mathrm{refl}}} \int_{ 0 }^{ T } V_{\alpha}^{\prime}(Z_{u}(t)) r(E) \, \mathrm{d}\ell_{Z_{u}}^{E}(t). 
\end{eqnarray} 
{\bf Martingale estimates}.
Define the sequence $M_{1}(\cdot )$ as 
\[
M_1(n) :=  \sum_{i=1}^{n}  \left( V_{\alpha}(Z_{u}(T_{u}(i))) - V_{\alpha}(E_{u}(i))  - P_{E_{u}(i)} \right), \quad n  \in \N_{>0} 
\] 
and the process $M_{2}(\cdot )$ as 
\[
M_2(T) := \int_{0}^{T} V^{\prime} _{\alpha}(Z_{u}(t))\sigma(Z_{u}(t),u(t))  \, \mathrm{d}W(t),  \quad T \ge 0.
\] 
In that notation, by using the boundary condition $V_{\alpha}^{\prime} (E) = 0$ for $E \in K_{\mathrm{refl}}$, we obtain 
\begin{equation} \label{e:JZ}
\int_{ 0 }^{T}  f(Z_{u}(t),u(t)) \, \mathrm{d}t \,- \sum_{E \in K_{\mathrm{res}} } P_{E}q_{u}(E, T) \le \alpha T - V_{\alpha}(Z_{u}(T)) + V_{\alpha}(0) + M_1(N_u(T)) + M_2(T).
\end{equation} 
The process $M_{1}(\cdot )$ is a martingale w.r.t. $\left( \mathcal{F}_{T_{u}(n)}  \right)_{n \in \N_{>0}} $, because the non-local boundary condition in  \eqref{e:HJB} implies
\begin{eqnarray*}
	&&\mathbb{E} \left(M_1(n+1) - M_1(n) \, | \, \mathcal{F}_{T_{u}(n)}\right)
      \\&& = \sum_{E \in K_{\mathrm{res}} }   \mathbb{P} \left(E_{u}(n+1) = E \, | \, \mathcal{F}_{T_{u}(n)}\right) \cdot \mathbb{E} \left(  V_{\alpha}(\gamma_{n+1}(E)) - V_{\alpha}(E) - P_{E}\, \Big| \, \mathcal{F}_{T_{u}(n)}  \right) 
      \\&& = \sum_{E \in K_{\mathrm{res}} }   \mathbb{P} \left(E_{u}(n+1) = E \, | \, \mathcal{F}_{T_{u}(n)}\right) \cdot \int_{\mathcal{X}} \left( V_{\alpha}(x) - V_{\alpha}(E) - P_{E} \right) \, \mathrm{d}\nu_{E}(x) = 0. 
\end{eqnarray*} 
By the strong law of large numbers for discrete martingales with uniformly bounded increments, it holds that
\begin{equation} \label{e:M1LN}
\lim_{n \to \infty} \frac{M_1(n)}{n} \stackrel{\text{a.s.}}{=} 0.
\end{equation} 
From Lemma \ref{l:T} and \eqref{e:M1LN}, it follows that
\begin{equation} \label{e:M1lil}
\lim_{T \to \infty} \frac{M_1(N_u(T))}{T}  \stackrel{\text{a.s.}}{=} 0.
\end{equation} 

Define the quadratic variation of the process $M_{2}(\cdot )$ as
\[
[M_{2}](T) := \int_{ 0 }^{ T }  \big( V^{\prime} _{\alpha}(Z_{u}(t)) \, \sigma(Z_{u}(t),u(t))  \big) ^2 \, \mathrm{d}t, \quad T \ge 0.	
\] 
We have $V^{\prime} _{\alpha} \in L_{\infty}[A,B]$ because $V_{\alpha} \in W^{2,1}[A,B]$ and $[A,B]$ is a bounded interval. 
According to Lemma \ref{l:equ}, there exists $H \in W^{2,1}[A,B]$ satisfying \eqref{e:H} with some $\Gamma$. By applying Itô decomposition for $H(Z_{u}(T))$ in the same way as in \eqref{e:Ito}, we obtain
\begin{equation} \label{e:Squa}
\int_{ 0 }^{T}  \sigma^2(Z_{u}(t),u(t)) \, \mathrm{d}t \, \le \Gamma T - H(Z_{u}(T)) + H(0) + \widetilde{M}_1(N_u(T)) + \widetilde{M} _2(T),
\end{equation} 
where 
\[
\widetilde{M}_1(n) :=  \sum_{i=1}^{n}  \big( H(Z_{u}(T_{u}(i))) - H(E_{u}(i))\big), \quad n  \in \N_{>0} 
\] 
and 
\[
\widetilde{M}_2(T) := \int_{0}^{T} H^{\prime}(Z_{u}(t))\sigma(Z_{u}(t),u(t))  \, \mathrm{d}W(t),  \quad T \ge 0.
\]
In a similar way as for $M_1(\cdot )$, using that the increments of $\widetilde{M}_{1}(\cdot ) $ are uniformly bounded, we can prove that
\begin{equation} \label{e:M1lln}
\limsup_{T \to \infty} \frac{\widetilde{M} _1(N_u(T))}{T}  \stackrel{\text{a.s.}}{\le} 0.
\end{equation} 
By the law of iterated logarithm for continuous local martingales, with probability 1, for sufficiently large $T$, we have 
\begin{equation} \label{e:M2lil}
	|\widetilde{M}_{2}(T)| \le 2\sqrt{ [\widetilde{M}_{2} ](T) \cdot  \log \log [\widetilde{M}_{2}](T)}, 
\end{equation} 
where the quadractic variation of the process $\widetilde{M}_{2} $ is
\begin{equation} \label{e:tM2}
[\widetilde{M}_{2}](T) := \int_{ 0 }^{ T }  \big( H^{\prime} (Z_{u}(t)) \, \sigma(Z_{u}(t),u(t))  \big) ^2 \, \mathrm{d}t \le  C_1 \int_{ 0 }^{ T }  \sigma^2(Z_{u}(t),u(t))\, \mathrm{d}t,
\end{equation}
for some $C_1 > 0$, due to $H^{\prime} \in AC[A,B]$.
By combining \eqref{e:Squa} together with $\eqref{e:M1lln}$, \eqref{e:M2lil}, and \eqref{e:tM2}, with probability 1, we obtain
\[
\frac{1}{C_1}[\widetilde{M}_{2}](T) \le \Gamma T + 2\sqrt{ [\widetilde{M}_{2} ](T) \cdot  \log \log [\widetilde{M}_{2}](T)} + O(T), \quad \text{as } T \to \infty.  
\]
Therefore, with probability 1, 
\[
[\widetilde{M}_{2}](T) = O(T),  \quad \text{as } T \to \infty.
\]
By combining this together with \eqref{e:Squa} and \eqref{e:M2lil}, we obtain, with probability 1, 
\[
\int_{ 0 }^{T}  \sigma^2(Z_{u}(t),u(t)) \, \mathrm{d}t = O(T),  \quad \text{as } T \to \infty.  
\] 
Therefore, the quadratic variation of $M_2(\cdot )$, with probability 1, also grows at most linearly
\[
[M_2](T) \le C_2 \int_{ 0 }^{ T }  \sigma^2(Z_{u}(t),u(t))\, \mathrm{d}t = O(T) ,  \quad \text{as } T \to \infty.
\] 
By the strong law of large numbers for continuous local martingales whose quadratic variation grows at most linearly, we have
\begin{equation} \label{e:M2LN}
\lim_{T \to \infty} \frac{M_2(T)}{T} \stackrel{\text{a.s.}}{=} 0.
\end{equation} 

By combining \eqref{e:JZ}, \eqref{e:M2LN}, and \eqref{e:M1lil} we obtain
\begin{equation*} \label{e:main}
J(u,f,P_{K_{\mathrm{res}}}) \stackrel{\text{a.s.}}{\le} \alpha
\end{equation*} 

{\bf Construction of the optimal control.}
By the measurable maximum theorem \cite[Th.~18.19]{AB06}, there exists a measurable selector $\eta^{*}(\cdot)$ satisfying \eqref{def:u0*}. 
Define
\[
	\mu^{*}(\cdot ) := \mu(\cdot , \eta^{*}(\cdot )), \quad
	\sigma^{*}(\cdot ) := \sigma(\cdot, \eta^{*}(\cdot )).
\]
By Condition \ref{cond:MSF_full}, $\mu^*(\cdot)$ and $\sigma^*(\cdot)$ satisfy the Engelbert--Schmidt condition \eqref{cond:ES}. Hence, for every initial condition $x\in \mathcal X$, equation \eqref{def:refl1} admits a weak solution with coefficients $\mu^*(\cdot)$ and $\sigma^*(\cdot)$ and reflecting boundaries in $K_{\mathrm{refl}}$.
Using these weak solutions together with the independent resetting variables $\gamma_n(E),E \in K_{\mathrm{res}}$, one can construct a filtered probability space, a Brownian motion, and a process $Z_{u^{*}}(\cdot )$ satisfying \eqref{def:Z_u} under the stationary Markov  control \eqref{def:u*}.

Under the control $u^{*}(\cdot )$, equality holds in \eqref{e:Ito_inequality} and \eqref{e:JZ} because of \eqref{def:u0*}. By applying \eqref{e:M2LN} and \eqref{e:M1lil}, we obtain \eqref{e:main_equality}.

\subsection{Proof of Theorem \ref{t:equ}}
First, we prove existence of solutions to the HJB equation \eqref{e:HJB} before imposing the boundary conditions for an arbitrary $\alpha \in \R$. Then, we use a shooting argument, to find $\alpha$ such that the solution to that auxiliary equation satisfies the non-local boundary condition. 

{\bf Solutions to the HJB equation before imposing the boundary condition}.
Fix $\alpha \in \R $. Consider the equation for a.e. $x  \in \left( A,B \right) $ 
\begin{equation} \label{e:HJB1}
	\begin{cases} 
		\alpha = \sup\limits_{\eta \in U } \left( \frac{1}{2}\sigma^2(x,\eta) V^{\prime \prime}_{\alpha}(x) + \mu(x,\eta) V^{\prime} _{\alpha}(x) + f(x,\eta)\right) , 
		\\ 
		V^{\prime} _{\alpha}(A)  = V_{\alpha}(A) = 0.   
	\end{cases}
\end{equation}
By the measurable maximum theorem \cite[Th.~18.19]{AB06}, there exists a measurable selector $u_1^{*}(\cdot ,\cdot ,\cdot )$ such that
\[
u_1^{*}(v^{\prime\prime},v^{\prime},x) \in  \argmax\limits_{\eta \in U} \left( \frac{1}{2}\sigma^2(x,\eta) v^{\prime \prime} + \mu(x,\eta) v^{\prime} + f(x,\eta)\right), \quad  (v^{\prime \prime},v^{\prime}, x )\in \R \times \R \times [A,B].
\] 
Evaluating the selector at $(v'',v',x) = (V_\alpha''(x),V_\alpha'(x),x)$ for a.e. $x \in (A,B)$, we obtain
\begin{equation*}
\alpha
=
\frac{1}{2}\sigma^2\bigl(x,u_1^*\bigr) V_\alpha''(x)
+
\mu\bigl(x,u_1^*\bigr) V_\alpha'(x)
+
f\bigl(x,u_1^*\bigr),
\end{equation*}
where $u_1^* := u_1^*\bigl(V_\alpha''(x),V_\alpha'(x),x\bigr).$ Since $\sigma^2(x,\eta) \ge \varsigma(x) > 0$ for all $(x,\eta) \in [A,B] \times U$, this identity determines $V_\alpha''(x)$ pointwise and yields the representation
\begin{equation*}
V_\alpha''(x) = g\bigl(V_\alpha'(x),x,\alpha\bigr),
\end{equation*}
for a measurable function $g:\R\times [A,B]\times \R \to \R$.
This function is dominated by
\begin{equation} \label{e:gL1}
|g(v^{\prime}, x,\alpha)| \le 2\left(\frac{|\alpha|}{\varsigma(x)} + \Lambda(x)|v^{\prime}|+F(x) \right) , \quad (v^{\prime}, x) \in \R \times [A,B] 
\end{equation}
and is Lipschitz in the first argument 
\begin{equation} \label{e:gL2}
|g(v^{\prime} _1,x,\alpha) - g(v^{\prime} _{2},x,\alpha)| \le 2\Lambda(x)|v^{\prime} _1-v^{\prime} _2|, \quad (v^{\prime}_{1},v^{\prime} _{2},x) \in \R \times \R \times  [A,B].
\end{equation} 
 By standard results on Carathéodory differential equations, from \eqref{e:gL1}, \eqref{e:gL2}, and Condition \ref{cond:MSF_full}, it follows that, for each $\alpha \in \R$, equation \eqref{e:HJB1} admits the unique solution $V_{\alpha} \in W^{2,1}[A,B]$. Moreover, the solution $V_{\alpha}(\cdot )$ depends continuously on the choice of $\alpha$, that is, if for $\left\{ \alpha_{i} \in \R \right\}_{i \in \N_{\ge 0}} $ it holds that $\alpha_{i} \stackrel{i \to \infty}{\longrightarrow} \alpha$, then $V_{\alpha_{i}} \rightrightarrows V_{\alpha}$ as $i \to \infty$.
Therefore, in order to find $\alpha \in \R$ for which $V_{\alpha}$ satisfies the boundary conditions of \eqref{e:HJB}, it suffices to show that the expression
\begin{equation} \label{cond:re}
	V_{\alpha}(B)+P_{B} - \int_{A}^{B} V_{\alpha}(x)  \, \mathrm{d}\nu_{B}(x)  
\end{equation}
takes both positive and negative values. 

{\bf The shooting argument}.
Define the function
\[
\Sigma(x) := \max_{\eta \in U}\sigma^2(x,\eta), \quad x \in [A,B].
\] 
Let $\varepsilon \in \left(0, \frac{B-A}{2}\right)$. Choose  $n_{\varepsilon} > 0$ such that, with 
\[
D_{\varepsilon} := \left\{ x \in [A,B] \, \Big| \,  2F(x)> n_{\varepsilon} \text{ or }  2\Lambda(x)> n_{\varepsilon}  \text{ or } \frac{2}{\Sigma(x)} < \frac{1}{n_{\varepsilon}}  \right\},
\] 
it holds that
\begin{eqnarray*}
	\int_{D_{\varepsilon}} 2 \Lambda(x) \, \mathrm{d}x < \frac{1}{2}, 
\qquad
	\int_{D_{\varepsilon}} 2F(x)  \, \mathrm{d}x < 1,
\end{eqnarray*} 
and
\begin{equation} \label{e:Deps}
|D_{\varepsilon}| < \varepsilon.
\end{equation} 
For a.e. $x \in (A,B)$, we have
\begin{equation}
V_{\alpha}^{\prime \prime}(x) \le 
\begin{cases} \label{e:V''}
2\left( \frac{\alpha}{\varsigma(x)}+ \Lambda(x)|V_{\alpha}^{\prime} (x)| + F(x) \right) , \quad \text{ if }  \frac{\alpha}{\varsigma(x)} + \Lambda(x)|V_{\alpha}^{\prime} (x)| + F(x) \ge 0, \\
\frac{2\varsigma(x)}{\Sigma(x)}\left(\frac{\alpha}{\varsigma(x)} + \Lambda(x)|V_{\alpha}^{\prime} (x)| + F(x)\right), \quad \text{ if }  \frac{\alpha}{\varsigma(x)} + \Lambda(x)|V_{\alpha}^{\prime} (x)| + F(x) < 0. 
\end{cases}  
\end{equation} 
 Let $k > 0$, $c_{k} := 2(k+1)$, $\alpha < -n_{\varepsilon}^2(1 + c_{k})$, $[x_0,x] \subset [A,B)$, and $y \in [x_0,x]$. If $y \notin D_{\varepsilon} $,  
then
\[
	\frac{2}{\varsigma(y)} \ge \frac{2}{\Sigma(y)} \ge \frac{1}{n_{\varepsilon}}.
\]
 If, in addition, $|V_{\alpha}^{\prime} (y)| \le c_{k}$, then  
\begin{equation} \label{e:nDe1}
	2\left( \frac{\alpha}{\varsigma(y)}+ \Lambda(y)|V^{\prime} _{\alpha}(y)|+F(y) \right)  \le -2n_{\varepsilon}(1+c_{k}) + n_{\varepsilon}c_{k} + n_{\varepsilon} < 0.
\end{equation}
This estimate together with \eqref{e:V''} implies that, for a.e. $y \notin D_{\varepsilon} $,
\begin{equation} \label{e:nDe2}
V^{\prime \prime}_{\alpha} (y) \le \frac{2\varsigma(y)}{\Sigma(y)}\left(\frac{\alpha}{\varsigma(y)} + \Lambda(y)|V^{\prime} _{\alpha}(y)|+F(y)\right) < 0.
\end{equation} 

 If $|V_{\alpha}^{\prime}(y)| \le c_{k} $ for all $y \in [x_0,x]$, then
\begin{equation} \label{e:V''+}
	\int_{x_0}^{x} V_{\alpha}^{\prime \prime}(y)   \, \mathrm{d}y \le \int_{[x_0,x]\cap D_{\varepsilon}} \big( 2F(y)+ 2\Lambda(y) |V^{\prime} _{\alpha}(y)| \big) \, \mathrm{d}y < 1 + \frac{c_{k}}{2} \le c_{k}-k. 
\end{equation} 
Therefore, $V_{\alpha}^{\prime}(x) < c_{k} $ for all $x \in [A,B]$, because $V_{\alpha}^{\prime} (A) = 0$. Moreover, if  $V_{\alpha}^{\prime}(y) > -c_{k} $ for all $y \in [x_0,x]$, then, due to \eqref{e:nDe1} and \eqref{e:nDe2}, it holds that
\begin{eqnarray} \nonumber
	V_{\alpha}^{\prime} (x) - V_{\alpha}^{\prime} (x_0) &<& c_{k} - k + \int_{[x_0,x] \setminus D_{\varepsilon}} 2\left( \frac{\alpha}{\varsigma(y)}+ \Lambda(y)|V^{\prime} _{\alpha}(y)|+F(y) \right) \, \mathrm{d}y 
	\\ &\le& \label{e:V'upper} 
	c_{k} - k + \left(\frac{\alpha}{n_{\varepsilon}} + n_{\varepsilon}c_{k} + n_{\varepsilon}\right) \big|[x_0,x] \setminus D_{\varepsilon}\big|.
\end{eqnarray} 
If $V^{\prime} _{\alpha}(x) = -c_{k}$ for some $x \in [A,B]$, then, by applying \eqref{e:V''+} on $[x,y]$ for some $y \in [x,B]$, we obtain 
\begin{equation} \label{e:V'upper2}
V^{\prime} _{\alpha}(y) < V_{\alpha}^{\prime}(x) + c_{k} - k < -k.
\end{equation} 
It follows from this, \eqref{e:Deps}, and \eqref{e:V'upper}, that for sufficiently large negative $\alpha$, we have 
\[
V_{\alpha}^{\prime} (x) < -k, \quad x \in [A+2\varepsilon, B].
\] 

Since $\nu_{B}$ is a probability distribution on $[A,B)$, there exists $\delta > 0$ such that 
\begin{equation} \label{e:dnu}
	\nu_{B}([A, B-\delta] ) > 0.
\end{equation}
Fix such a $\delta$.
Choose $k>0$ sufficiently large and then choose
$\varepsilon>0$ sufficiently small so that $B-\delta>A+2\varepsilon$ and
\[
(k\delta-2c_k\varepsilon)\nu_B([A,B-\delta])>P_B.
\]
For sufficiently large negative $\alpha$, we have
$V_\alpha'(x)<-k$ for $x\in[A+2\varepsilon,B]$. Hence, for
$x\in[A,B-\delta]$,
\[
V_\alpha(x)-V_\alpha(B)\ge k\delta-2c_k\varepsilon.
\]
Therefore
\[
\int_A^B V_\alpha(x)\,d\nu_B(x)-V_\alpha(B)
\ge
(k\delta-2c_k\varepsilon)\nu_B([A,B-\delta])
>P_B.
\]
It follows that, for sufficiently large negative $\alpha$, the expression
\eqref{cond:re} takes negative values. The case of sufficiently large positive $\alpha$ is obtained by analogous (symmetric) estimates and yields positive values of \eqref{cond:re}. This concludes the proof of the theorem.

\subsection{Proof of Lemma \ref{l:equ}}
In the case of one resetting and one reflecting boundary, Theorem \ref{t:equ} implies that the HJB equation \eqref{e:HJB} with $\sigma^2(\cdot )$ instead of $f(\cdot )$ and $P_{B} = 0$ admits a solution  $V_{\alpha} \in W^{2,1}[A,B]$ for some $\alpha \in \R$. Theorem \ref{t:equ} can be applied because Condition
\ref{cond:MSF_full} remains valid with $\sigma^2(\cdot )$ in place of $f(\cdot )$. Therefore, equation \eqref{e:H} holds with $\Gamma = \alpha$ and $H = V_{\alpha}$. 

Let $K_{\mathrm{res}} = \{A,B\}$. For $\Gamma\in \R$ consider a solution $H_{\Gamma} \in W^{2,1}[A,B]$ to the equation, for a.e. $x  \in ( A,B ) $, 
\begin{equation*} \label{e:HJB2}
	\begin{cases} 
		\Gamma = \sup\limits_{\eta \in U } \left( \frac{1}{2}\sigma^2(x,\eta) H_{\Gamma}^{\prime \prime}(x) + \mu(x,\eta) H_{\Gamma}^{\prime}(x) + \sigma^2(x,\eta)\right), 
		\\ 
		H_{\Gamma}^{\prime}(A)  = -4, \quad H_{\Gamma}(A)=0.   
	\end{cases}
\end{equation*}
Existence follows by the same argument as in the proof of Theorem \ref{t:equ}. It is enough to find $\Gamma \in \R$ such that $H_{\Gamma}(B) = 0$ and $H_{\Gamma}(x) \le 0$ for all $x \in (A,B)$.  
Let $x_{\circ}(\Gamma) := \inf\left\{ x \in (A,B] \, | \, H_{\Gamma}(x) = 0 \right\} $.
 $x_{\circ}(\Gamma)$ is decreasing in $\Gamma$, because, if $\Gamma_1 > \Gamma_{2}$, then $H_{\Gamma_1} \ge H_{\Gamma_2}$ on $[A,B]$.
Using the shooting argument as in the proof of Theorem \ref{t:equ} with $k = 1$, from \eqref{e:V'upper} and \eqref{e:V'upper2}, we obtain that, for sufficiently large negative $\Gamma$, 
\[
H^{\prime} _{\Gamma}(x) \le -1, \quad x \in [A,B].
\] 
Therefore, for sufficiently large negative $\Gamma$, we have $H_{\Gamma}(x) \le 0 $ for all $x \in [A,B]$. Using an analogous argument, we obtain that, for sufficiently large $\Gamma$, we have $H_{\Gamma}(B) > 0$.
Since the family of solutions $\left( H_{\Gamma} \right)_{\Gamma \in \R} $ depends continuously on the parameter $\Gamma$, the map $\Gamma\mapsto H_\Gamma(B)$ is continuous. By the monotonicity of
$x_\circ(\Gamma)$ and the preceding estimates, the value of $\Gamma$ can be
chosen so that $x_\circ(\Gamma)=B$.

\subsection{Proof of Lemma \ref{l:T}} \label{s:l:T}
Define functions $\Delta, \widetilde{\Delta}: [A,B] \to \R$ by equations
\begin{equation} \label{e:Dinf}
\inf_{\eta \in U } \left( \frac{1}{2} \sigma^2(x,\eta) \Delta^{\prime \prime}(x)  + \mu(x,\eta) \Delta^{\prime}(x)    \right) = -1, 
\end{equation} 
\begin{equation} \label{e:Dsup}
\sup_{\eta \in U } \left( \frac{1}{2} \sigma^2(x,\eta) \widetilde{\Delta} ^{\prime \prime}(x)  + \mu(x,\eta) \widetilde{\Delta} ^{\prime}(x)    \right) = -1, 
\end{equation} 
with boundary conditions
\begin{eqnarray*}
	&&	\Delta(E) = \widetilde{\Delta}(E) = 0, \quad E \in K_{\mathrm{res}}, 
	\\
	&& \Delta^{\prime} (E) = \widetilde{\Delta}^{\prime} (E) = 0, \quad E \in K_{\mathrm{refl}}.  
\end{eqnarray*}
The existence of $\Delta, \widetilde{\Delta} \in W^{2,1}[A,B]$ follows by the same one-dimensional shooting argument as in the proof of Theorem \ref{t:equ}, applied to equations \eqref{e:Dinf} and \eqref{e:Dsup} with the corresponding boundary conditions. 
Since for $x \in \mathcal{X}$, \ $\Delta(x)$ is the expected lifetime for the extremal diffusion associated with the infimum in \eqref{e:Dinf} started from $x$ before reaching $K_{\mathrm{res}}$, we have $\Delta(x) >0$.
Define the times between consecutive boundary hits as
\[
	\tau_{u}(n) := T_{u}(n) - T_{u}(n-1), \quad n \in \N_{>0}
\]
and, for $\tau \ge 0$, define a random variable 
\[
	\Delta_{u}(n,\tau) := \mathbb{E} \left(\left( \tau_{u}(n+1) - \tau \right)_{+} \, \Big| \, \mathcal{F}_{ T_{u}(n) + \tau} \right),
\]
where $y_{+} := \max (y,0)$ for $y \in \R$.
From \eqref{e:Dinf}, \eqref{e:Dsup}, on the event $\left\{ \tau \le \tau_{u}(n+1) \right\}$, Dynkin's formula applied on the shifted interval $[T_{u}(n)+\tau, T_{u}(n+1)]$ yields
\begin{eqnarray*}
	 \label{e:dnu} 
	 \widetilde{\Delta}(Z_{u}(T_{u}(n) + \tau)) \ge \Delta_{u}(n,\tau) \ge \Delta(Z_{u}(T_{u}(n) + \tau)), \quad \text{a.s.} 
 \end{eqnarray*}
From this and that $\widetilde{\Delta} \in L_{\infty}[A,B]$, there exists $\widetilde{c}  > 0$ such that
\begin{equation} \label{e:D_upper}
\widetilde{c} \ge \Delta_{u}(n,\tau), \quad \tau \in [0, \tau_{u}(n+1)], \quad \text{a.s.} 
\end{equation} 
Since $\left\{ \nu_{E} \, | \, E \in K_{\mathrm{res}} \right\} $ are probability distributions on $\mathcal{X}$ and $\Delta(\cdot ) > 0$ on $\mathcal{X} $, there exist $\varepsilon> 0$ 
such that the conditional expectations of $\tau_{u}(n)$ are uniformly separated from 0 
\begin{eqnarray} \nonumber
	\mathbb{E} \left(\tau_{u}(n+1)\, | \, \mathcal{F}_{T_{u}(n)-} \right) &=&
	\mathbb{E} \left(\Delta_{u}(n,0) \, | \, \mathcal{F}_{T_{u}(n)-}\right) 
	\\ &\ge& \nonumber
	 \mathbb{E} \big(\Delta(Z_{u}(T_{u}(n))) \, \big| \, \mathcal{F}_{T_{u}(n)-}\big)
	\\ &=& \nonumber 
	\int_{\mathcal{X}} \Delta(x) \, \mathrm{d}\nu_{E_u(n)}(x)
	\\ &\ge& \label{e:eps^2}\min_{E \in K_{\mathrm{res}}} \int_{\mathcal{X}} \Delta(x) \, \mathrm{d}\nu_{E}(x) \ge \varepsilon.
\end{eqnarray} 
To prove that the reset times grow at least linearly almost surely, define a sequence 
\[
	M_3(n) := \sum_{i=1}^{n}   \left(\tau_{u}(i+1) - \mathbb{E} \left(\tau_{u}(i+1)  \mid \mathcal{F}_{T_{u}(i)-}\right)  \right), \quad n \in \N_{>0}. 
\] 
The process $M_3(\cdot )$ is a martingale w.r.t. filtration $\left\{ \mathcal{F}_{T_{u}(n+1)-} \right\}_{n \in \N_{> 0}} $ such that
\[
	\mathbb{E} \left(M_3(n)-M_3(n-1)\right)^2 < 2\widetilde{c} ^{\, \, 2}, \quad n \in \N_{>1}
\]
because, by using \eqref{e:D_upper}, we obtain,
\begin{eqnarray*}
	\mathbb{E} \left(   \tau_{u}(n+1)^2 \, \big| \, \mathcal{F}_{T_{u}(n)-}  \right) &=& \mathbb{E}   \left( \int_{ 0 }^{\tau_{u}(n+1)}  2\left( \tau_{u}(n+1) - t \right)   \, \mathrm{d}t \, \Big| \,  \mathcal{F}_{T_{u}(n)-}\right)   
										   \\ &=& \mathbb{E}  \left( \int_{ 0 }^{\tau_{u}(n+1)}   2 \, \mathbb{E} \left( \tau_{u}(n+1) - t \, \Big| \, \mathcal{F}_{T_{u}(n)+t}\right)   \, \mathrm{d}t \, \Big| \, \mathcal{F}_{T_{u}(n)-} \right)  	\\ &=& \mathbb{E}  \left( \int_{ 0 }^{\tau_{u}(n+1)}   2 \Delta_{u}\left(n, t\right)    \, \mathrm{d}t \, \Big| \, \mathcal{F}_{T_{u}(n)-} \right)  
		     \\ &\le& 2\widetilde{c} ^{\, \, 2}.
\end{eqnarray*} 
It follows from this and the strong law of large numbers for martingales whose predictable quadratic variation grows at most linearly that 
\begin{equation} \label{e:M3lil}
\lim_{n \to \infty} \frac{M_3(n)}{n} \stackrel{\text{a.s.}}{=} 0.
\end{equation} 
By combining \eqref{e:eps^2} and \eqref{e:M3lil}, with probability 1, we obtain that, for sufficiently large $n \in \N_{>0}$,
\[
T_{u}(n) = \sum_{i=1}^{n} \tau_{u}(i) \ge \frac{\varepsilon}{2} n.  
\] 

\section{Proofs of the results for renewal diffusions} \label{s:prd}

\subsection{Proof of Theorem \ref{t:B} and Lemma \ref{l:d}} 
Define the random variables
\[
	Y_{i}(f) := \int_{T_{i-1}}^{ T_{i} } f(Z(t))  \, \mathrm{d}t, \quad Y_{i}(1) := T_{i} - T_{i-1}, \quad i \in \N_{> 0}.
\]
The random variables $\left( Y_{i}(f) \right)_{i \in \N_{>1}}$ are independent and identically distributed, as are the random variables $\left( Y_{i}(1) \right)_{i \in \N_{>1}}$. The random variables $Y_1(f)$ and $Y_{1}(1)$ may have a different distribution from those of the remaining sequence, but it is negligible in the long-run averages.

Lemma \ref{l:d} follows from the strong law of large numbers for $\left\{ Y_{i}(1) \right\}_{i \in \N_{\ge 0}} $, since $Y_{i}(1) > 0$ a.s. 
Since either $\mathbb{E}  Y_2(f)<\infty$ or $\mathbb{E}  Y_2(1)<\infty$, the strong law of large numbers implies that 
\[
	\lim_{n \to \infty} \frac{\sum_{i=1}^{n} Y_{i}(f)}{\sum_{i=1}^{n} Y_{i}(1) } \stackrel{\text{a.s.}}{=} \frac{\mathbb{E} Y_{2}(f)}{\mathbb{E} Y_{2}(1)}, \quad \lim_{n \to \infty} \frac{n}{\sum_{i=1}^{n} Y_{i}(1) } \stackrel{\text{a.s.}}{=} \frac{1}{\mathbb{E} Y_{2}(1)}.   
\] 
Since $S(B) < \infty$, it follows from \cite[Ch.~II.11]{BS02} and \cite[p.~160]{IM74} that, for all $x \in (A,B)$, it holds that 
\[
\mathbb{E} \left( \int_{ T_{1} }^{ T_2 } f(Z(t)) \, \mathrm{d}t \, \Big| \, Z(T_1) = x \right)  = G_{B}(x,f).
\] 
It follows that
\[
\mathbb{E} Y_{2}(f) = G_{B}(f), \quad \mathbb{E} Y_{2}(1) = G_{B}(1),
\] 
possibly infinite.
From this and the fact that $\sum_{i=1}^{n} Y_{i}(1) = T_{n}$, we obtain 
\begin{equation} \label{e:fLLN0}
\lim_{n \to \infty} \frac{\int_{ 0 }^{T_{n}} f(Z(t))  \, \mathrm{d}t}{T_{n}} \stackrel{\text{a.s.}}{=} \frac{G_{B}(f)}{G_{B}(1)}, \quad \lim_{n \to \infty} \frac{n}{T_{n}} \stackrel{\text{a.s.}}{=} \frac{1}{G_{B}(1)}.
\end{equation} 
For $n \in \N_{>0}$ and any $T \in [T_{n}, T_{n+1})$, with probability 1, we have
\begin{equation} \label{e:f+-}
\frac{\int_{ 0 }^{T_{n}} f(Z(t))  \, \mathrm{d}t}{T_{n+1}} \le \frac{\int_{ 0 }^{T} f(Z(t))  \, \mathrm{d}t }{T} \le \frac{\int_{ 0 }^{T_{n+1}} f(Z(t))  \, \mathrm{d}t}{T_{n}}
\end{equation} 
 If $G_{B}(1) < \infty$, then from \eqref{e:fLLN0} it follows that the left-hand side and the right-hand side of \eqref{e:f+-}, with probability 1, converge to the same value (possibly infinity). If $G_{B}(1) = \infty$ and $G_{B}(f) < \infty$, then both sides converge to 0. Indeed,
 \[
\frac{\int_{ 0 }^{T_{n+1}} f(Z(t))  \, \mathrm{d}t}{T_{n}} = \frac{\frac{1}{n+1}\int_{ 0 }^{T_{n+1}} f(Z(t))  \, \mathrm{d}t}{\frac{1}{n+1}T_{n}} \stackrel{\text{a.s.}}{\to} 0. 
 \]
 The lower bound is handled in the same way.
 In both cases, combining this with the second limit in \eqref{e:fLLN0} yields the desired result
\begin{equation*} \label{e:JfX}
J(f,P_{B} ) \stackrel{\text{a.s.}}{=} \frac{G_{B}(f)-P_{B}}{G_{B}(1)}.
\end{equation*}

\begin{rem}
When $A>-\infty$, Theorem \ref{t:B} can be proved using Theorem \ref{t:main} in the same way as Theorem \ref{t:AB}.
\end{rem}
\subsection{Proof of Theorem \ref{t:AB}}
Since $A > -\infty$ and $|S(E)| < \infty$ for $ E \in \left\{ A,B \right\} $, we have
\[
G_{AB}^{E}( 1) < \infty, \quad E \in \left\{ A,B \right\} .
\] 
Since $\frac{f}{\sigma^2} \in L_{1}[A,B]$, it follows that
\[
G_{AB}^{E}(f ) < \infty, \quad E \in \left\{ A,B \right\}.
\] 
With $\alpha$ defined by \eqref{e:JfAB}, the function  
\[
V(x) := V(A)\frac{S(B)-S(x)}{S(B)-S(A)} + V(B)\frac{S(x)-S(A)}{S(B)-S(A)} + G_{AB}(x,f-\alpha), \quad x \in [A,B] 
\] 
is a $W^{2,1}[A,B]$ solution to \eqref{e:HJB} prior to imposing the boundary conditions. Choosing $V(A) := 0$ and
\[
V(B) := \frac{P_{A} - G_{AB}^{A}(f) + \alpha\,G_{AB}^{A}(1)}{p_{AB}}, 
\] 
ensures that the non-local boundary conditions in \eqref{e:HJB} are satisfied.
Therefore, Theorem \ref{t:main} yields \eqref{e:JfAB}.
\subsection{Proof of Theorem \ref{t:special}}

If condition \eqref{cond:fs} holds, then, due to \cite[Th.~2.6]{UM12} and the regularity of $X$, with probability 1, for all sufficiently large $n \in \N_{>0}$, we have
\[
\int_{ 0 }^{T_n} f(Z(t))   \, \mathrm{d}t = \infty.
\] 
If condition \eqref{cond:fS} holds, then, due to \cite[Th.~2.11]{UM12}, we have
\[
\int_{ 0 }^{T_1} f(Z(t))   \, \mathrm{d}t \stackrel{\text{a.s.}}{=} \infty.
\] 

Suppose that $G_B(1) = \infty$. This occurs only when $A = -\infty$.
It holds that 
\[
\lim_{T \to \infty} \frac{-P_{B} \cdot N(T)}{T} \stackrel{\text{a.s.}}{=} \lim_{n \to \infty} -\frac{P_{B}\cdot n}{T_{n}} \stackrel{\text{a.s.}}{=}  -\frac{P_{B}}{\mathbb{E} T_1} = -\frac{P_{B}}{G_{B}(1)} = 0.
\]
Therefore, whenever it exists, $J(\cdot ,P_{B} )$ does not depend on $P_{B}$, namely, for any $h: (-\infty, B) \to \R_{\ge 0}$, 
\begin{equation} \label{e:ind}
	J(h,P_{B} ) \stackrel{\text{a.s.}}{=}  \lim_{T \to \infty} \frac{1}{T}\int_{ 0 }^{ T } h(Z(t))  \, \mathrm{d}t .
\end{equation} 
Assume that neither condition \eqref{cond:fs} nor \eqref{cond:fS} holds, and that $f_{-\infty} := \lim_{x \to -\infty} f(x)$ exists.  Let  $a < B$ and define
\[
	1_{a}(x) := 1_{(-\infty, a]}(x), \qquad f_{a}(x) := f(x)1_{a}(x), \quad x \in (-\infty,B).
\]
Then $G_{B}(f-f_{a}) < \infty$ and $G_{B}(1-1_{a}) < \infty$. From this and $G_{B}(1) = \infty$, by Theorem \ref{t:B}, we obtain 
\begin{equation} \label{e:rest}
J(f-f_{a},P_{B}) \stackrel{\text{a.s.}}{=} 0, \quad J(1-1_{a}, P_{B})\stackrel{\text{a.s.}}{=} 0.
\end{equation}
From \eqref{e:ind}, \eqref{e:rest}, and $\int_{ 0 }^{ T } f = \int_{ 0 }^{ T } f_{a} + \int_{ 0 }^{ T } (f-f_{a})$ for $T \ge 0$,
we obtain
\[
J(f_{a},P_{B}) \stackrel{\text{a.s.}}{=} J(f,P_{B}), \quad J(1_{a}, P_{B})\stackrel{\text{a.s.}}{=} 1.
\] 
Since $f(x) \to f_{-\infty}$ as $ x \to -\infty$, it holds that  
\begin{equation} \label{e:limsup}
\limsup_{a \to -\infty} \limsup_{T \to \infty} \frac{1}{T} \int_{ 0 }^{ T } f_{a}(Z(t))  \, \mathrm{d}t \le f_{-\infty}, \quad \text{a.s.}
\end{equation} 
and
\begin{equation} \label{e:liminf1}
	\liminf_{a \to -\infty} \liminf_{T \to \infty} \frac{1}{T} \int_{ 0 }^{ T } f_{a}(Z(t))  \, \mathrm{d}t \ge f_{-\infty} \cdot \liminf_{a \to -\infty} J(1_{a},P_{B}) \ge f_{-\infty}, \quad \text{a.s.} 
\end{equation} 
From \eqref{e:limsup} and \eqref{e:liminf1}, it follows that
\[
J(f,P_{B} ) \stackrel{\text{a.s.}}{=} f_{-\infty}.
\] 
\section{Proofs of the optimization examples} \label{s:pfoe}
\subsection{Proof of Example \ref{exa:1}}
Fix $\mu > 0$. The scale function and the density of the speed measure of the Brownian motion with a drift $\mu$ are given by
\[
S_{\mu}(x) = \frac{1 - e^{-2\mu x}}{2\mu}, \quad m_{\mu}(x) = e^{2\mu x}, \quad x \in \R. 
\] 
By substituting these expressions into \eqref{def:GKE}, we obtain
\[
	G_{B}(\mu, 1) = \frac{B}{\mu}, 
\] 
and
\begin{equation} \label{e:GBmuf}
G_{B}(\mu,f) = \int_{0}^{B}\frac{e^{-2\mu x} - e^{-2\mu B}}{\mu B} \, xe^{2\mu x}\, \mathrm{d} x = \frac{2\mu^2 B^2 - 2\mu B + 1 - e^{-2\mu B}}{4\mu^3 B}. 
\end{equation} 
By Theorem \ref {t:B}, it holds that
\[
\mathbb{E}  J(\mu,f,P_{B}) = \frac{2\mu^2 B^2 - 2\mu B + 1 - e^{-2\mu B}}{4\mu^2 B^2}- \frac{P_{B}\mu}{B}.
\] 
We maximize this expression over $\mu > 0$. We have 
\[
\frac{\mathrm{d}}{\mathrm{d}\mu} \mathbb{E}  J(\mu,f,P_{B}) = \frac{B\mu(1+e^{-2\mu B}) + e^{-2\mu B} - 1}{2\mu^3 B^2} - \frac{P_{B}}{B}.
\] 
Therefore, the optimal $\mu$ satisfies 
\begin{equation} \label{e:mu_0}
\frac{B\mu(1+e^{-2\mu B}) + e^{-2\mu B} - 1}{2\mu^3 B} = P_{B}.
\end{equation} 
The left-hand side of that equality is decreasing on $(0, \infty)$ and has limit $\frac{B^2}{3}$ as $\mu \to 0$. Hence, for $P_{B} = 0$, the objective functional $\mathbb{E} J(\mu,f,0)$ 
is increasing in $\mu$ and equality \eqref{case:0} holds. 
For $P_{B} \in (0, \frac{B^2}{3})$, the equation \eqref{e:mu_0} has a unique solution, which maximizes $\mathbb{E} J(\mu,f, P_{B})$. The equation \eqref{case:2} is obtained by introducing the P\'eclet number $\Pe := \mu B$.
For $P_{B} \ge \frac{B^2}{3}$, the objective functional $\mathbb{E} J(\mu,f,P_{B})$ is negative for all $\mu > 0$ and equals 0 for $\mu = 0$. The latter follows from $G_B(0,1) = \infty$ and $G_{B}(0,f) < \infty$.

\subsection{Proof of Example \ref{exa:2}}
For $\mu \in \R$, we have
\[
G_{B}(\mu,1) = 2\int_0^{B} \frac{1 - e^{-2\mu(B-x)}}{2\mu} \, \mathrm{d}x = 
\frac{2\mu B-1+e^{-2\mu B}}{2\mu^2}, 
\] 
and $G_{B}(\mu,f)$ has the same form as in \eqref{e:GBmuf}.
Then, according to Theorem \ref {t:B}, it holds that 
\[
\mathbb{E} J(\mu,f,P_{B}) = \frac{G_{B}(\mu, f) - P_{B}}{G_{B}(\mu,1)} = 
\frac{(2\Pe^2 - 2\Pe + 1)e^{2\Pe} - 1}{2\Pe((2\Pe - 1)e^{2\Pe} + 1)} - \frac{P_{B}}{B^2} \cdot  
\frac{2\Pe^2}{2\Pe - 1 + e^{-2\Pe}},
\] 
where $\Pe := \mu B$.

\section{Proof of Theorem \ref{t:infinite}} \label{s:pi}
A trivial upper bound is
\begin{equation} \label{est:I_up}
      I(\mathcal{D}_{AB},f) \le \esssup f.
\end{equation} 

In order to get the opposite estimate, we will find a diffusion which, by the ergodic theorem, spends most of the time at the points where
$f(\cdot )$ takes large values.

Consider a diffusion $\widetilde{X}(\cdot ) $ with state space $(A, B)$, some continuous drift coefficient $\mu(\cdot )$ that will be defined later, and the constant diffusion coefficient $\sigma \equiv 1$. Its speed measure has the following simplified form 
\begin{equation} \label{def:m}
    m_{\widetilde{X} }(x) = \exp\left(2 \int_{0}^{x} \mu(s) \, \mathrm{d}s\right), \quad x \in (A,B),
\end{equation}
and the form of the scale function is 
\[
   S_{\widetilde{X} }(x) := \int^{x}_{0} \frac{1}{m_{\widetilde{X} }(y)}\, \mathrm{d}y \quad x \in (A,B).
\] 
Assume that $-\infty < S_{\widetilde{X} }(x) < \infty$ for all $x \in (A,B)$.
The levels $A$ and $B$ belong to the {\it entrance boundary} in the Feller classification and do not belong to the {\it exit boundary} under condition
\begin{equation} \label{cond:border}
    S_{\widetilde{X} }(B) = \int^{B}_{0} \frac{1}{m_{\widetilde{X} }(x)} \, \mathrm{d}x 
    = -S_{\widetilde{X} }(A) 
    = \int_{A}^{0} \frac{1}{m_{\widetilde{X} }(x)}\, \mathrm{d}x = \infty.
\end{equation} 
A diffusion $\widetilde{X}(\cdot ) $ is called {\it recurrent} if, for all $x,y \in (A, B)$, it is true that
\[
    \mathbb{P} \left(\min \left\{ t \, | \, \widetilde{X}_{x}(t) = y \right\} < \infty \right) = 1,
\] 
where $\widetilde{X}_{x}$ is a diffusion that has the same drift and diffusion coefficients as $\widetilde{X}(\cdot )$, but is started at $x$ rather than at 0. According to \cite[Ch.~IV.12, Prop.~12.1]{aB17}, diffusion $\widetilde{X}(\cdot ) $ is recurrent if condition \eqref{cond:border} holds.  

Assuming that the speed measure of the diffusion $\widetilde{X}(\cdot ) $ is finite, denote 
\begin{equation} \label{def:p}
     p_{\widetilde{X}}(x) := Q^{-1}m_{\widetilde{X} }(x), \quad x \in (A,B),
\end{equation}
where $Q := \int_{A}^{B} m_{\widetilde{X}} (x) \mathrm{d}x$.
According to \cite[Ch.~IV.11, p.~319]{aB17}, the distribution with density $p_{\widetilde{X}}(\cdot )$ is the stationary distribution of the diffusion $\widetilde{X} (\cdot )$.

For arbitrary $n>0$, define a function 
$f_{n}:(A, B) \to \R_{\ge 0}$ as 
\[
     f_{n}(x) := \min (f(x), n), \quad x \in (A, B).
\] 
The finiteness of $Q$ together with the recurrence of $\widetilde{X}(\cdot )$ are the sufficient conditions, under which 
the ergodic theorem for one-dimensional diffusions holds, see \cite[Ch.~V.7.53]{RW00} and \cite[Ch.~6.8]{IM74}. In particular,
\begin{equation} \label{e:ergo}
	J(\widetilde{X},f_{n} ) = \mathbb{E} \left( \lim_{T \to \infty} \frac{1}{T} \int_{0}^{T}f_{n}(\widetilde{X} (t)) \, \mathrm{d}t \right)  = \int_{A}^{B} f_{n}(x) p_{\widetilde{X}}(x) \, \mathrm{d}x.
\end{equation}
We will now specify the drift of a quasi-optimal diffusion.
Since $f_{n} \in L_{\infty}(A, B)$ and is non-negative, we have
\[
   \sup \left\{ \int_{A}^{B} f_{n}(x)g(x) \, \mathrm{d}x \, \Big| \, {g: (A,B) \to \R_{\ge 0}, \lVert g \rVert_{L_1} \le 1} \right\}
   = \lVert f_{n} \rVert_{L_{\infty}} 
   = \esssup f_{n}. 
\] 
Let $\varepsilon > 0$. Since continuously differentiable functions are dense in the unit $L_{1}$-ball \\
$\left\{ h \in L_1(A,B) \, | \, \lVert h \rVert_{L_1} \le 1 \right\} $,
there exists a function of this class $g_{\varepsilon}(\cdot )$ satisfying 
$\lVert g_{\varepsilon} \rVert_{L_1} \le 1$ and 
\begin{equation} \label{e:delta}
    \int_{A}^{B} f_{n}(x)g_{\varepsilon}(x) \, \mathrm{d}x \ge \esssup f_{n}- \varepsilon.
\end{equation}
Moreover, we may assume, without loss of generality, that $g_{\varepsilon}(\cdot )$ is {\it strictly}
positive.
Let $\delta > 0$ be such that
\[
     \int_{A+\delta}^{B-\delta}  f_{n}(x)  g_{\varepsilon}(x)  \, \mathrm{d}x 
     \ge \esssup f_{n}- 2\varepsilon.
\] 
Modify the function $g_{\varepsilon}(\cdot )$ to a continuously differentiable, strictly positive function $g_{\varepsilon,1}(\cdot )$ so that 
\[
\begin{cases}
g_{\varepsilon,1}(x) = g_{\varepsilon}(x), \quad x \in (A+\delta, B-\delta),
\\
g_{\varepsilon,1}(x) \le g_{\varepsilon}(x), \quad x \in (A, A+\delta) \cup (B-\delta, B),
\end{cases}
\] 
and condition \eqref{cond:border} holds with $g_{\varepsilon,1}(\cdot )$ in place of $m_{\widetilde{X} }(\cdot )$, i.e.
\begin{equation} \label{cond:border2}
      \int^{B}_{0} \frac{1}{ g_{\varepsilon,1}(x)} \, \mathrm{d}x  
    = \int_{A}^{0} \frac{1}{ g_{\varepsilon,1}(x)}\, \mathrm{d}x = \infty.
\end{equation} 

Finally, define the drift coefficient of a diffusion $\widetilde{X}(\cdot ) $ as 
\begin{equation} \label{def:mu}
	\mu(x) := \frac{g_{\varepsilon,1}^{\,\prime}(x) }{2g_{\varepsilon,1} (x)}, \quad x \in (A,B). 
\end{equation}
Under this choice,  by  \eqref{def:mu},\eqref{def:m}, and \eqref{def:p},   the functions $g_{\varepsilon,1}(\cdot )$, $m_{\widetilde{X}}(\cdot )$, and $p_{\widetilde{X}}(\cdot )$ coincide up to some multiplicative constants. In particular, we have 
$p_{\widetilde{X}}(\cdot ) = \tfrac{g_{\varepsilon,1}(\cdot )}{||g_{\varepsilon,1}||_{L_1}}$.

Using $||g_{\varepsilon,1}||_{L_1} \le ||g_{\varepsilon}||_{L_1}\le 1$, as well as \eqref{e:ergo} and \eqref{e:delta}, we obtain
\begin{eqnarray*}
    I(\mathcal{D}_{AB},f) &\ge&  I(\mathcal{D}_{AB},f_{n}) \ge J(\widetilde{X},f_{n}) 
    = \int_{A}^{B} f_{n}(x)  \frac{g_{\varepsilon,1}}{||g_{\varepsilon,1}||_{L_1}} \, \mathrm{d}x 
\\
    &\ge&  \int_{A}^{B} f_{n}(x) g_{\varepsilon,1}(x) \, \mathrm{d}x
    \ge  \int_{A+\delta}^{B-\delta} f_{n}(x) g_{\varepsilon,1}(x) \, \mathrm{d}x
    =  \int_{A+\delta}^{B-\delta} f_{n}(x) g_{\varepsilon}(x) \, \mathrm{d}x
\\    
   &\ge& \esssup f_{n} - 2\varepsilon.
\end{eqnarray*} 
Letting $\varepsilon \downarrow 0$ and $n \to \infty$ successively, we obtain
\begin{equation} \label{est:I_low}
I(\mathcal{D}_{AB},f) \ge \esssup f.
\end{equation} 
We conclude by combining estimates \eqref{est:I_up} and \eqref{est:I_low}
\[
     I(\mathcal{D}_{AB},f) = \esssup f.
\]

\bigskip
The author is grateful to M.A. Lifshits for setting the problem and for valuable assistance with this work.

\end{document}